\documentclass[12pt]{amsproc}%{amsart}
\usepackage{amsmath,amssymb,amscd}

\begin{document}

\title[Free Probability for Pairs of Faces II]{Free Probability for Pairs of Faces II: \\
$2$-Variables Bi-free Partial $R$-Transform and Systems with Rank $\le 1$ Commutation}
\author{Dan-Virgil Voiculescu}
\address{D.V. Voiculescu \\ Department of Mathematics \\ University of California at Berkeley \\ Berkeley, CA\ \ 94720-3840}
\thanks{Research supported in part by NSF Grant DMS-1301727.}
\date{}

\begin{abstract}
We compute the generating series for the simplest class of bi-free cumulants, beyond the free cumulants, the two-bands bi-free cumulants of a pair of left and a right variable. We also consider two-faced systems with a commutation condition implying that two-bands moments, that is expectation values of the product of a monomial of left and a monomial of right variables determine the other moments. Examples include hyponormal operators, dual systems in free entropy theory and bi-partite systems.
\end{abstract}

\maketitle

\setcounter{section}{-1}
\section{Introduction}
\label{sec0}

In Part~I (\cite{10}) we introduced bi-freeness, which underlies the extension of free probability to two-faced systems of non-commutative random variables, that is, to systems with ``left'' and ``right'' variables. We proved on the basis of general Lie-theory considerations the existence of bi-free cumulants and found the analogue of Gaussian variables and distributions in this theory.

Here in Part~II our focus will be on the bi-freeness machinery to handle a wider class of distributions of left and right variables, which had been an important motivation for us to study bi-freeness and which goes beyond bi-free Gaussians. These are distributions where algebraic relations between left and right variables insure that all moments can be computed from the ``two-bands'' moments $\varphi(LR)$ where $L$ and $R$ are monomials in left and, respectively, right variables. The systems satisfy commutation relations between left and right variables a feature which is preserved under bi-free addition. The simplest example is provided by the bi-partite systems where left and right variables commute with each other. Other examples like those arising from dual systems in free entropy theory (\cite{9}) and from hyponormal operators (\cite{1},\cite{3}) are characterized by commutators of rank $\le 1$ between left and right variables. To get a hold on additive bi-free convolution for such systems it suffices to compute their ``two-bands'' bi-free cumulants. In the simplest case of one left variable and one right variable we find a formula connecting the two-variables generating series of two-bands moments and the partial $R$-transform which is the two-variables generating series of the two-bands cumulants. This solves the problem of computing bi-free additive convolutions of the distributions of hermitian and antihermitian parts of normal operators or more generally of hyponormal operators with rank $\le 1$ self-commutator. It may be interesting to note that two-bands moments from pairs of variables arising from normal or hyponormal operators appear also in other questions in mathematical physics (\cite{4}).

Not counting the introduction, this paper has four sections. Section~\ref{sec1} contains preliminaries about bi-freeness and bi-free cumulants from Part~I and about notation we will use. Section~\ref{sec2} gives our main result about the partial $R$-transform for a pair consisting of one left variable and one right variable. This is analogous to our result for the usual one-variable $R$-transform \cite{8}. Since neither our original approach via a canonical Toeplitz operator form for the variable (\cite{8}) nor a combinatorial machinery based on partitions (\cite{5}) were at hand, it was very convenient that one of the proofs of Haagerup (\cite{2}) could be used as a starting point for the bi-free result. Then in Section~\ref{sec3} we discuss systems of left and right variables satisfying a rank $\le 1$ commutation condition, their distributions and behavior under bi-free addition.

Section~\ref{sec4} deals with the characterization of classical independence of left and right variables by the vanishing of the bi-free cumulants which mix left and right indices.

\section{Preliminaries}
\label{sec1}

\noindent
{\bf 1.1. Bi-freeness.} Let $({\mathcal A},\varphi)$ be a non-commutative probability space and $((a_i^{(k)})_{i \in I},(b_j^{(k)})_{j \in J})$, $k = 1,2$ be two-faced families of non-commutative random variables in $({\mathcal A},\varphi)$. Then the definition in $2.6.$~Part~I of bi-freeness of the two-faced families for $k = 1$ and $k = 2$ amounts in view of the results in \S2 of Part~I to the existence of two vector spaces with specified state vector $({\mathcal X}_k,\overset{\circ}{{\mathcal X}}_k,\xi_k)$, $k = 1,2$ with the following properties. There are ${\tilde a}_i^{(k)} \in {\mathcal L}({\mathcal X}_k)$, ${\tilde b}_j^{(k)} \in {\mathcal L}({\mathcal X}_k)$, for $i \in I$, $j \in J$, $k = 1,2$ so that if $({\mathcal X},\overset{\circ}{{\mathcal X}},\xi) = ({\mathcal X}_1,\overset{\circ}{{\mathcal X}}_1,\xi_1) * ({\mathcal X}_2,\overset{\circ}{{\mathcal X}}_2,\xi_2)$ and $\lambda_k$, $\rho_k$ are the left and right representations of ${\mathcal L}({\mathcal X}_k)$ on ${\mathcal L}({\mathcal X})$, then the joint distribution of the variables $a_i^{(k)}$, $b_j^{(k)}$, $i \in I$, $j \in J$, $k \in \{1,2\}$ in $({\mathcal A},\varphi)$ is the same as the joint distribution in $({\mathcal L}({\mathcal X}),\varphi_{\xi})$ of the variables $\lambda_k({\tilde a}_i^{(k)})$, $\rho_k({\tilde b}_j^{(k)})$, $i \in I$, $j \in J$, $k = \{1,2\}$. Moreover, one can always choose ${\mathcal X}_k = {\mathcal A}$, $\overset{\circ}{{\mathcal X}}_k = \ker \varphi$, $\xi_k = 1$, $k = 1,2$ and ${\tilde a}_i^{(k)}$, ${\tilde b}_j^{(k)}$ to be the left multiplication operators on ${\mathcal A}$ by ${\tilde a}_i^{(k)}$ and, respectively, ${\tilde b}_j^{(k)}$. This is due to the fact that in the definition any $({\mathcal X}_k,\overset{\circ}{{\mathcal X}}_k,\xi_k)$ and ${\tilde a}_i^{(k)},{\tilde b}_j^{(k)}$ yield the same joint distribution of the family as long as the joint distribution of each $({\tilde a}_i^{(k)})_{i \in I} \cup ({\tilde b}_j^{(k)})_{j \in J}$ in $({\mathcal L}({\mathcal X}_k),\varphi_{\xi_k})$ equals that of $(a_i^{(k)})_{i \in I} \cup (b_j^{(k)})_{j \in J}$ in $({\mathcal A},\varphi)$.

\bigskip
\noindent
{\bf 1.2. Bi-free cumulants.} If $I$ and $J$ are index sets for left and, respectively, right variables of a two-faced family of non-commutative random variables, we obtained in Theorem~5.7 of Part~I a general result about existence and uniqueness of bi-free cumulants $R_{\alpha}$ indexed by maps $\alpha: \{1,\dots,n\} \to I \coprod J$, $n \ge 1$. Here $R_{\alpha}$ is a homogeneous polynomial of degree $n$ in commutative variables $X_{\alpha K}$ or degree $|K|$, where the index $\alpha K$ denotes the map $\alpha \square b_K: \{1,\dots,|K|\} \to I \coprod J$, where $\emptyset \ne K \subset \{1,\dots,n\}$ is a subset and $b_K: \{1,\dots,|K|\} \to K$ the increasing bijection. Equivalently the $R_{\alpha}$'s are indexed by sequences $\alpha = (h_1,h_2,\dots,h_n) \in (I \coprod J)^n$ and are homogeneous polynomials of degree $n$ in the commuting variables $X_{h_{k_1},\dots,h_{k_r}}$ of degree~$r$, where $(h_{k_1},\dots,h_{k_r})$, $1 \le k_1 < \dots < k_r \le n$ runs over the subsequences of $h_1,\dots,h_n$. If $z = ((z_i)_{i \in I},(z_j)_{j \in J})$ is a two-faced family in $({\mathcal A},\varphi)$, then $R_{\alpha}(z)$ or $R_{\alpha}(\mu_z)$ where $\mu_z$ is the distribution of $z$, will denote the value of $R_{\alpha}$ when we give the $X_{\alpha(k_1),\dots,\alpha(k_r)}$ the values $\varphi(z_{\alpha(k_1)}\dots z_{\alpha(k_r)})$ where $1 \le k_1 < \dots < k_r \le n$. The polynomial $R_{\alpha}$ is characterized by the requirements that the coefficient of the leading term in $X_{\alpha(1),\dots,\alpha(n)}$ equals $1$ and that the cumulant additivity property
\[
R_{\alpha}(z'+z'') = R_{\alpha}(z') + R_{\alpha}(z'')
\]
is satisfied if $z'$ and $z''$ are bi-free.

In the present paper we will be interested in perhaps the simplest bi-free cumulants which are not free cumulants. These are those where $n = p+q$, $p > 0$, $q > 0$ and there are $\beta: \{1,\dots,p\} \to I$, $\gamma: \{1,\dots,q\} \to J$ so that $\alpha: \{1,\dots,p+q\} \to I \coprod J$ is given by $\alpha(k) = \beta(k)$ if $1 \le k \le p$ and $\alpha(p+l) = \gamma(l)$ if $1 \le l \le q$. We shall denote $R_{\alpha}$ by $R_{\beta\gamma}$ in this particular case.

Actually, we will be studying in detail the case of $R_{\beta\gamma}$, when $|I| = |J| = 1$. In this situation it will be sufficient to specify $p$ and $q$ and we will use the notation $R_{pq}$. This is the case of one left and one right variable $a$ and $b$. Then $R_{pq}$ is a polynomial in variables $X_{rs}$, $0 \le r \le p$, $0 \le s \le q$, $r+s > 0$ and $R_{pq}(a,b)$ corresponds to setting $X_{rs}$ equal to $\varphi(a^rb^s)$. Here the coefficient of $X_{pq}$ is $1$ and $R_{pq}$ is completely determined by the additivity property for the addition of bi-free pairs.

\section{The partial bi-free $R$-transform $R_{(a,b)}(z,w)$ for a pair of variables}
\label{sec2}

\noindent
{\bf 2.1.} Let $(a,b)$ in $({\mathcal A},\varphi)$ be a two-faced pair of non-commutative random variables, that is $a$ will be a left variable and $b$ will be a right variable. We will prove a formula for the generating series
\[
R_{(a,b)}(z,w) = \sum_{\substack{
m \ge 0 \\
n \ge 0 \\
m + n \ge 1
}} R_{m,n}(a,b)z^mw^n
\]
where $R_{m,n}(a,b)$ are bi-free cumulants for $(a,b)$ of the type discussed in the section of preliminaries. In particular, this means $R_{m,n}(a,b)$ will be a polynomial in moments $\varphi(a^pb^q)$ where $0 \le p \le m$, $0 \le q \le n$, $p + q \ge 1$. We will show how to derive $R_{(a,b)}(z,w)$ from a two-variables Green's function
\[
\varphi((z1-a)^{-1}(w1-b)^{-1}) = G_{(a,b)}(z,w).
\]
The result will also involve the one-variable Green's functions
\[
\varphi((z1-a)^{-1}) = G_a(z),\ \varphi((z1-b)^{-1}) = G_b(z)
\]
and the formula for one-variable $R$-transforms $R_a(z) = K_a(z) - z^{-1}$, $R_b(z) = K_b(z) - z^{-1}$ we found in \cite{8}, where $K_a(z),K_b(z)$ are the inverses of $G_a(z),G_b(z)$ at $\infty$. In establishing a formula for $R_{(a,b)}(z,w)$ we will use in an essential way insights from one of the proofs given by U.~Haagerup \cite{2} for our formula for the one-variable $R$-transform. It may be interesting to note that certain two-variable Green's functions play an important role in the theory of hyponormal operators (\cite{1},\cite{3}).

Throughout Section~\ref{sec2} $({\mathcal M},\varphi)$ will denote a non-commutative probability space.

\bigskip
\noindent
{\bf 2.2. Lemma.} {\em Let $(({\mathcal A}_i,{\mathcal B}_i))_{i \in I}$ be a bi-free family of included pairs of faces in $({\mathcal M},\varphi)$. Let $a_k \in {\mathcal A}_{\alpha(k)}$, $b_l \in {\mathcal B}_{\beta(l)}$, $1 \le k \le m$, $1 \le l \le n$ be so that $\varphi(a_k) = \varphi(b_l) = 0$ and $\alpha(1) \ne \alpha(2) \ne \dots \ne \alpha(m)$, $\beta(1) \ne \beta(1) \ne \dots \ne \beta(n)$. Then we have
\[
\varphi(a_m\dots a_1b_n \dots b_1) = \delta_{m,n} \prod_{1 \le k \le \min(m,n)} \delta_{\alpha(k),\beta(k)}\varphi(a_kb_k).
\]
}

\bigskip
\noindent
{\bf {\em Proof.}} Let $({\mathcal X}_i,\overset{\circ}{{\mathcal X}}_i,\xi_i)$, $i \in I$ be vector spaces with specified state vector and let $\mu_i: {\mathcal A}_i \to {\mathcal L}({\mathcal X}_i)$, $\nu_i: {\mathcal B}_i \to {\mathcal L}({\mathcal X})$ be homomorphisms so that
\[
\begin{split}
&\varphi(a_m\dots a_1b_n \dots b_1) = p(\lambda_{\alpha(m)}(\mu_{\alpha(m)}(a_m)) \\
&\dots \lambda_{\alpha(1)}(\mu_{\alpha(1)}(a_1)) \\
&\dots \rho_{\beta(n)}(\nu_{\beta(n)}(b_n)) \\
&\dots \rho_{\beta(1)}(\nu_{\beta(1)}(b_1)))
\end{split}
\]
as required by the bi-freeness assumption.

We shall proceed by induction over $\min(m,n)$. The case of $\min(m,n) = 0$ and $\delta_{m,n} = 0$ is the case of $m > 0 = n$ or $m = 0 < n$, that is there are either only left or only right variables in the moment and the moment is then $0$ by the definition of freeness. The case of $m = 0 = n$ is the case of the empty word, which should be viewed as being $\varphi(1) = 1$.

Assume now $\min(m,n) > 0$. Then
\[
\begin{split}
&\rho_{\beta(n)}(\nu_{\beta(n)}(b_n)) \dots \rho_{\beta(1)}(\nu_{\beta(1)}(b_1))\xi \\
&= \nu_{\beta(1)}(b_1)\xi_{\beta(1)} \otimes \dots \otimes \nu_{\beta(n)}(b_n)\xi_{\beta(n)} \in \overset{\circ}{{\mathcal X}}_{\beta(1)} \otimes \dots \otimes \overset{\circ}{{\mathcal X}}_{\beta(n)}.
\end{split}
\]
If $\alpha(1) \ne \beta(1)$, then
\[
\begin{split}
&\lambda_{\alpha(m)}(\mu_{\alpha(m)}(a_m))\dots \lambda_{\alpha(1)}(\mu_{\alpha(1)}(a_1))\rho_{\beta(n)}(\nu_{\beta(n)}(b_n)) \dots \rho_{\beta(1)}(\nu_{\beta(1)}(b_1))\xi \\
&= \mu_{\alpha(m)}(a_m)\xi_{\alpha(m)} \otimes \dots \otimes \mu_{\alpha(1)}(a_1)\xi_{\alpha(1)} \otimes \\
&\qquad \nu_{\beta(1)}(b_1)\xi_{\beta(1)} \otimes \dots \otimes \nu_{\beta(n)}(b_n)\xi_{\beta(n)} \\
&\in \overset{\circ}{{\mathcal X}}_{\alpha(m)} \otimes \dots \otimes \overset{\circ}{{\mathcal X}}_{\alpha(1)} \otimes \overset{\circ}{{\mathcal X}}_{\beta(1)} \otimes \dots \otimes \overset{\circ}{{\mathcal X}}_{\beta(n)}
\end{split}
\]
and hence $\varphi(a_m\dots a_1b_n \dots b_1) = 0$. This proves the assertion in case $\alpha(1) \ne \beta(1)$, that is $\delta_{\alpha(1),\beta(1)} = 0$.

If $\alpha(1) = \beta(1)$, then
\[
\begin{split}
&\lambda_{\alpha(m)}(\mu_{\alpha(m)}(a_m)) \dots \lambda_{\alpha(1)}(\mu_{\alpha(1)}(a_1))(\nu_{\beta(1)}(b_1)\xi_{\beta(1)} \otimes \dots \otimes \nu_{\beta(n)}(b_n)\xi_{\beta(n)} \\
&= \lambda_{\alpha(m)}(\mu_{\alpha(m)}(a_m))\dots \lambda_{\alpha(2)}(\mu_{\alpha(2)}(a_2))(\mu_{\alpha(1)}(a_1)\nu_{\beta(1)}(b_1)\xi_{\beta(1)} - \varphi(a_1b_1)\xi_{\beta(1)}) \\
&\qquad\otimes \nu_{\beta(2)}(b_2)\xi_{\beta(2)} \otimes \dots \otimes \nu_{\beta(n)}\xi_{\beta(n)}) \\
&+ \varphi(a_1b_1)\lambda_{\alpha(m)}(\mu_{\alpha(m)}(a_m)\dots \lambda_{\alpha(2)}(\mu_{\alpha(2)}(a_2))\rho_{\beta(n)}(\nu_{\beta(n)}(b_n)) \dots \rho_{\beta(1)}(\nu_{\beta(1)}(b_1)) \\
&\qquad \in \overset{\circ}{{\mathcal X}}_{\alpha(m)} \otimes \dots \otimes \overset{\circ}{{\mathcal X}}_{\alpha(2)} \otimes \overset{\circ}{{\mathcal X}}_{\beta(1)} \otimes \dots \otimes \overset{\circ}{{\mathcal X}}_{\beta(n)} \\
&+ \varphi(a_1b_1)\lambda_{\alpha(m)}(\mu_{\alpha(m)}(a_m)) \dots \lambda_{\alpha(2)}(\mu_{\alpha(2)}(a_2))\rho_{\beta(n)}(\nu_{\beta(n)}(b_n))\dots\rho_{\beta(2)}(\nu_{\beta(2)}(b_2))\xi
\end{split}
\]
and hence
\[
\varphi(a_m\dots a_1b_n \dots b_1) = \varphi(a_1b_1)\varphi(a_m\dots a_2b_n \dots b_2).
\]
The assertion then follows from the induction hypothesis applied to
$\varphi(a_m \dots a_2b_n \dots b_2)$. \qed

\bigskip
\noindent
{\bf 2.3.} We introduce in this subsection some of the assumptions and notation we shall use to derive the main result of the section. To make use of some of the technique in \cite{2}, we will bring as far as possible, our notations close to those in \cite{2}. We will, however, avoid the precise inequalities for radii of convergence of power-series, by using the language of germs of holomorphic functions. These will be often germs of holomorphic functions near $0$ or $\infty$ taking values in ${\mathbb C}$ or in some Banach space.

We will assume $({\mathcal A},\varphi)$ our non-commutative probability space is a unital Banach algebra ${\mathcal A}$ and that the expectation functional $\varphi$ is continuous. We will deal with two two-faced pairs of variables $(a_1,b_1)$ and $(a_2,b_2)$ in $({\mathcal A},\varphi)$ which will be added (note that here our notation differs from \cite{2} where the variables to be added are $a$ and $b$).

We define $h_{a_k}(t_k) = \varphi((1-t_ka_k)^{-1})$, $h_{b_k}(s_k) = \varphi((1-s_kb_k)^{-1})$, $k = 1,2$, are germs of holomorphic functions near $0$ with values in ${\mathbb C}$. Further we shall also consider $a_k(t_k) = (1-t_ka_k)^{-1} - h_{a_k}(t_k)1$, $b_k(s_k) = (1 - s_kb_k)^{-1} - h_{b_k}(s_k)1$, $k = 1,2$, which are germs of holomorphic functions near $0$ taking values in ${\mathcal A}$.

We shall also use some of our usual notation (\cite{8},\cite{11}) for the one-variable $R$-transform. So $G_{a_k}(z_k) = \varphi((z_k1-a_k)^{-1})$, $G_{b_k}(z_k) = \varphi((z_k1-b_k)^{-1})$ will be germs of holomorphic functions near $\infty$ while $G_a(t^{-1}) = \varphi((t^{-1}1-a)^{-1})$, etc., will be germs of holomorphic functions near $0$. Then $K_{a_k}(z_k)$ will be a germ of meromorphic function near $0$, taking the value $\infty$ at $0$, which is the inverse of $G_{a_k}(z_k)$, etc. The $R$-transform $R_{a_k}(z_k) = K_{a_k}(z_k) - z_k^{-1}$ is then the germ of a holomorphic function near $0$.

We will also need two-variable analogs of $G_a(z)$, namely the
\[
G_{(a,b)}(z,w) = \varphi((z1-a)^{-1}(w1-b)^{-1})
\]
which is a germ of holomorphic function near $(\infty,\infty)$ in $({\mathbb C} \cup \{\infty\}) \times ({\mathbb C} \cup \{\infty\})$. We will also use the related
\[
H_{(a,b)}(t,s) = \varphi((1-ta)^{-1}(1-sb)^{-1})
\]
which is a germ of holomorphic function near $(0,0) \in {\mathbb C} \times {\mathbb C}$.

\bigskip
\noindent
{\bf 2.4. Theorem.} {\em We have the following equality of germs of holomorphic functions near $(0,0) \in {\mathbb C}^2$
\[
R_{(a,b)}(z,w) = 1 + zR_a(z)+wR_b(w) - \frac {zw}{G_{a,b}(K_a(z),K_b(w))}
\]
if $(a,b)$ is a two-faced pair in a Banach-algebra non-commutative probability space.
}

\bigskip
\noindent
{\bf 2.5.} The proof of the theorem will occupy the rest of Section~\ref{sec2} and will consist of several steps to establish two facts. The first will be, that if $(a_1,b_1)$ and $(a_2,b_2)$ are bi-free in $({\mathcal A},\varphi)$ then
\[
\begin{split}
&\left( \frac {zw}{G_{a_1,b_1}(K_{a_1}(z),K_{b_1}(w))} - 1\right) \\
&+ \left( \frac {zw}{G_{a_2,b_2}(K_{a_2}(z),K_{b_2}(w))} - 1\right) \\
&= \frac {zw}{G_{a_1+a_2,b_1+b_2}{(K_{a_1+a_2}(z),K_{b_1+b_2}(w))}} - 1
\end{split}
\]
as germs of holomorphic functions near $(0,0)$. The second will be, that the coefficient of $z^mw^n$, $m \ge 0$, $n \ge 0$, $m+n > 0$ of the series expansion of
\[
1 + 2R_a(z) + wR_b(w) - \frac {zw}{G_{a,b}(K_a(z),K_b(w))}
\]
is a polynomial of bi-degree $(m,n)$ in the $\varphi(a^kb^l)$, where $\varphi(a^kb^l)$ is given bi-degree $(k,l)$, and that the coefficient of the leading term $\varphi(a^mb^n)$ is $1$. Clearly the additivity part will be the more difficult part of the proof, the checking of the bi-degree properties being routine.

\bigskip
\noindent
{\bf 2.6.} To compute $G_{(a_1+a_2,b_1+b_2)}(z,w)$ when $(a_1,b_1)$ and $(a_2,b_2)$ are bi-free we shall first compute $H_{(a_1+a_2,b_1+b_2)}(t,s)$. Here $(z,w)$ will be near $(\infty,\infty)$, while $(t,s)$ will be near $(0,0)$. An argument in \cite{2} can then be used to transform $H_{(a_1+a_2,b_1+b_2)}(t,s) = \varphi((1-t(a_1+a_2))^{-1}(1-s(b_1+b_2))^{-1})$ into an expectation value of the form
\[
\begin{split}
&\varphi((1-t_2a_2)^{-1}(1-\rho a_1(t_1)a_2(t_2))^{-1}(1-t_1a_1)^{-1}(1-s_2b_2)^{-1} \\
&(1-\sigma b_1(s_1)b_2(s_2))^{-1}(1-s_1b_1)^{-1})
\end{split}
\]
for $t_1,t_2,s_1,s_2$ close to $0$ and $\rho,\sigma$ sufficiently close to $1$. The last quantity can be further transformed expanding $(1-\rho a_1(t_1)a_2(t_2))^{-1}$ and $(1-\sigma b_1(s_1)b_2(s_2))^{-1}$ in geometric series and expressing the $(1-ta)^{-1}$, $(1-sb)^{-1}$ as $a(t)+h_a(t)1$ and $b(s)+h_b(s)1$. Thus we get to a computation which will be recorded as the next lemma.

\bigskip
\noindent
{\bf 2.7. Lemma.} {\em Assuming $(a_1,b_1)$ and $(a_2,b_2)$ are bi-free, we have
\[
\begin{split}
&\varphi((a_2(t_2)+h_{a_2}(t_2)1)(a_1(t_1)a_2(t_2))^k(a_1(t_1)+h_{a_1}(t_1)1)(b_2(s_2)+h_{b_2}(s_2)1) \\
&\qquad (b_1(s_1)b_2(s_2))^l(b_1(s_1)+h_{b_1}(s_1)1)) \\
&=\delta_{k,l}(\varphi(a_1(t_1)b_1(s_1)))^k(\varphi(a_2(t_2)b_2(s_2)))^l(\varphi(a_1(t_1)b_1(s_1)) \\
&+ h_{a_1}(t_1)h_{b_1}(s_1))(\varphi(a_2(t_2)b_2(s_2))+h_{a_2}(t_2)h_{b_2}(s_2))
\end{split}
\]
for $s_1,s_2,t_1,t_2$ near $0$ and $k \ge 0$, $l \ge 0$.
}

\bigskip
\noindent
{\bf {\em Proof.}} We start by noticing that among the 16 terms which arise from expanding the left-hand side, only the following four may be non-zero.
\[
\begin{split}
&\varphi((a_2(t_2)a_1(t_1))^{k+1}(b_2(s_2)b_1(s_1))^{l+1}) \\
&+ \varphi(a_2(t_2)(a_1(t_1)a_2(t_2))^kb_2(s_2)(b_1(s_1)b_2(s_2))^l)h_{a_1}(t_1)h_{b_1}(s_1) \\
&+ \varphi((a_1(t_1)a_2(t_2))^ka_1(t_1)(b_1(s_1)b_2(s_2))^lb_1(s_1))h_{a_2}(t_2)h_{b_2}(s_2) \\
&+ \varphi((a_1(t_1)a_2(t_2))^k(b_1(s_1)b_2(s_2))^l)h_{a_1}(t_1)h_{a_2}(t_2)h_{b_2}(s_2).
\end{split}
\]
Indeed since $\varphi(a_j(t_j)) = \varphi(b_i(s_i)) = 0$ and the expectation values are such that Lemma~2.2 applies. First, the number of left and right factors, that is of $a_j(t_j)$'s and $b_i(s_i)$'s must be equal, so this accounts for the vanishing of 8 terms in the expansion in which there is an even number of factors of one kind and an odd number of the other kind. Then there are four more terms which vanish because the sequence of indices 1 and 2 for the $a$'s must be the same as for the $b$'s. The application of Lemma~2.2 then gives that the four terms, which may not vanish, are actually equal to
\[
\begin{split}
&\delta_{k,l}((\varphi(a_2(t_2)b_2(s_2))^{k+1}(\varphi(a_1(t_1)b_1(s_1)))^{k+1} \\
&+ (\varphi(a_2(t_2)b_2(s_2)))^{k+1}\varphi(a_1(t_1)b_1(s_1))^k)h_{a_1}(t_1)h_{b_1}(s_1) \\
&+ (\varphi(a_1(t_1)b_1(s_1)))^{k+1}(\varphi(a_2(t_2)b_2(s_2)))^kh_{a_2}(t_2)h_{b_2}(s_2) \\
&+ (\varphi(a_1(t_1)b_1(s_1))^k(\varphi(a_2(t_2)b_2(s_2)))^kh_{a_1}(t_1)h_{a_2}(t_2)h_{b_1}(s_1)h_{b_2}(s_2),
\end{split}
\]
which is just the right-hand side of the equality we wanted to prove.\qed

\bigskip
\noindent
{\bf 2.8.} We will use the previous result in conjunction with expansions of $(1-\rho a_1(t_1)a_2(t_2))^{-1}$ and $(1-\sigma b_1(s_1)b_2(s_2))^{-1}$. The reader should keep in mind for this that as $t$ and $s$ tend to zero also $a_k(t),b_k(s)$ tend to zero. Thus, if $\rho$ and $\sigma$ are bounded, the expansions converge for $t_k,s_k$ near zero.

\bigskip
\noindent
{\bf 2.9. Lemma.} {\em Assume $(a_1,b_1)$ and $(a_2,b_2)$ are bi-free and $|\rho| \le C$, $|\sigma| \le C$ for some constant $C > 0$. Then we have
\[
\begin{split}
&\varphi(((1-t_1a_1)(1-\rho a_1(t_1)a_2(t_2))(1-t_2a_2))^{-1} \\
&\qquad ((1-s_1b_1)(1-\sigma b_1(s_1)b_2(s_2))(1-s_2b_2))^{-1}) \\
&= \frac {(\varphi(a_1(t_1)b_1(s_1)+h_{a_1}(t_1)h_{b_1}(s_1))(\varphi(a_2(t_2)b_2(s_2))+h_{a_2}(t_2)h_{b_2}(s_2))}{1 - \rho\sigma\varphi(a_1(t_1)b_1(s_1))\varphi(a_2(t_2)b_2(s_2))}
\end{split}
\]
for $t_1,t_2,s_1,s_2$ near $0$.
}

\bigskip
\noindent
{\bf {\em Proof.}} The left-hand side equals
\[
\begin{split}
&\varphi((a_2(t_2)+h_{a_2}(t_2)1)(1-\rho a_1(t_1)a_2(t_2))^{-1}(a_1(t_1)+h_{a_1}(t_1)1) \\
&\qquad(b_2(s_2)+h_{b_2}(s_1)1)(1-\sigma b_1(s_1)b_2(s_2))^{-1}(b_1(s_1)+h_{b_1}(s_1)1))
\end{split}
\]
which in view of 2.8 and Lemma~2.7 is further equal to
\[
\begin{split}
&\sum_{k \ge 0,l \ge 0} \rho^k\sigma^l\varphi((a_2(t_2)+h_{a_2}(t_2)1)(a_1(t_1)a_2(t_2))^k(a_1(t_1)+h_{a_1}(t_1)1) \\
&\qquad (b_2(s_1)+h_{b_2}(s_2)1)(b_1(s_1)b_2(s_2))^l(b_1(s_1)+h_{b_2}(s_1)1)) \\
&= \sum_{k \ge 0} \rho^k\sigma^k(\varphi(a_1(t_1)b_1(s_1)))^k(\varphi(a_2(t_2)b_2(s_2))^k \\
&\cdot (\varphi(a_1(t_1)b_1(s_1)+h_{a_1}(t_1)h_{b_1}(s_1))(\varphi(a_2(t_2)b_2(s_2))+h_{a_2}(t_2)h_{b_2}(s_2))).
\end{split}
\]
The lemma then follows from
\[
\begin{split}
&\sum_{k\ge 0} \rho^k\sigma^k(\varphi(a_1(t_1)b_1(s_1)))^k(\varphi(a_2(t_2)b_2(s_2)))^k \\
&= (1-\rho\sigma\varphi(a_1(t_1)b_1(s_1))\varphi(a_2(t_2)b_2(s_2)))^{-1}.
\end{split}
\]
\qed

\bigskip
\noindent
{\bf 2.10.} The next Lemma is a result from \cite{2} (part of Lemma~3.3 and Theorem~3.4, see pages 139--141 in \cite{2}). For the reader's convenience a proof will be sketched.

\bigskip
\noindent
{\bf 2.11. Lemma.} {\em {\em a)} If $t_1,t_2$ are near zero and $t_1h_{a_1}(t_1) = t_2h_{a_2}(t_2)$, then $\rho = (h_{a_1}(t_1)h_{a_2}(t_2))^{-1}$,
\[
t = \frac {t_1h_{a_1}(t_1)}{h_{a_1}(t_1)+h_{a_2}(t_2)-1}
\]
are such that
\[
(1-t_1a_1)(1-\rho a_1(t_1)a_2(t_2))(1-t_2a_2) = \frac {h_{a_1}(t_1)+h_{a_2}(t_2)-1}{h_{a_1}(t_1)h_{a_2}(t_2)} (1-t(a_1+a_2)).
\]

{\em b)} If $t$ sufficiently close to zero is given then $t_1,t_2$ sufficiently close to zero exist so that $t_1h_{a_1}(t_1) = t_2h_{a_2}(t_2)$ and
\[
t = \frac {t_1h_{a_1}(t_1)}{h_{a_1}(t_1)+h_{a_2}(t_2)-1}
\]
and $|\rho| < 2$ where $\rho = (h_{a_1}(t_1)h_{a_2}(t_2))^{-1}$.
}

\bigskip
\noindent
{\bf {\em Proof.}} a) We have
\[
\begin{split}
&(1-t_1a_1)(1-\rho a_1(t_1)(a_2(t_2))(1-t_2a_2) \\
&= (1-t_1a_1)(1-t_2a_2) - \rho(1-h_{a_1}(t_1)1+t_1h_{a_1}(t_1)a_1) \\
&\qquad(1-h_{a_2}(t_2)1+t_2h_{a_2}(t_2)a_2) \\
&= c_0 + c_1a_1 + c_2a_2 + c_3a_1a_2,
\end{split}
\]
where
\[
\begin{split}
c_3 &= t_1t_2(1-\rho h_{a_1}(t_1)h_{a_2}(t_2)) = 0 \\
c_0 &= 1-\rho(h_{a_1}(t_1)-1)(h_{a_2}(t_2)-1) = \frac {h_{a_1}(t_1)+h_{a_2}(t_2)-1}{h_{a_1}(t_1)h_{a_2}(t_2)} \\
c_1 &= -t_1(1+\rho h_{a_1}(t_1)-\rho h_{a_1}(t_1)h_{a_2}(t_2)) = -t_1(h_{a_2}(t_2))^{-1} \\
c_2 &= -t_2(1+\rho h_{a_2}(t_2)-\rho h_{a_1}(t_1)h_{a_2}(t_2)) = -t_2(h_{a_1}(t_1))^{-1}.
\end{split}
\]
In view of the assumption that $t_1h_{a_1}(t_1) = t_2h_{a_2}(t_2)$, we have that $c_1 = c_2$ and, moreover,
\[
t = \frac {t_1h_{a_1}(t_1)}{h_{a_1}(t_1)+h_{a_2}(t_2)-1} = c_0^{-1} \cdot t_1(h_{a_2}(t_2))^{-1} = -c_1c_0^{-1} = -c_2c_0^{-1}.
\]
Thus we have
\[
\begin{split}
&(1-t_1a_1)(1-\rho a_1(t_1)a_2(t_2))(1-t_2a_2) \\
&= c_0(1 - ta_1 - ta_2) = \frac {h_{a_1}(t_1)+h_{a_2}(t_2)-1}{h_{a_2}(t_2)h_{a_2}(t_2)} (1-t(a_1+a_2)).
\end{split}
\]

b) Since $t_k \to t_kh_{a_k}(t_k)$ is a holomorphic map near zero taking value $0$ at $0$ and with first order derivative equal $1$ at $0$, there is $p(z)$ holomorphic near $0$, $p(0) = 0$, $p'(0) = 1$ so that $t_1h_{a_1}(t_1) = p(t_1)h_{a_2}(p(t_1))$ for $t_1$ near zero. Thus, given $t$ close to zero we must show there is $t_1$ so that
\[
t = \frac {t_1h_{a_1}(t_1)}{h_{a_1}(t_1)+\frac {t_1}{p(t_1)}h_{a_1}(t_1)-1}.
\]
The right-hand side of this equation is a holomorphic function of $t_1$ near zero, with first-order jet $(0,1)$, hence there is a solution $t_1$ which is a holomorphic function of $t$ near zero with first-order jet $(0,1)$ at $0$. Also taking $t_2 = p(t_1)$, $\rho$ becomes a holomorphic function of $t_1$ near zero, which takes the value $1$ at $0$. In particular, we have $|\rho| < 2$ if $t_1$ is close enough to $0$.\qed

\bigskip
\noindent
{\bf 2.12.} Clearly a statement similar to the previous lemma also holds with $b_k,s_k,s$ replacing $a_k,t_k,t$. The same remark should also be made about the next lemma, which provides further details from \cite{2} to the previous lemma.

\bigskip
\noindent
{\bf 2.13. Lemma.} {\em If $a_1,a_2$ are free in $({\mathcal A},\varphi)$ and $t,t_1,t_2$ are near zero and satisfy the relations in Lemma~$2.11$, then we have: $h_{a_1+a_2}(t) = h_{a_1}(t_1)+h_{a_2}(t_2)-1$ and $th_{a_1+a_2}(t) = t_1h_{a_1}(t_1) = t_2h_{a_2}(t_2)$.
}

\bigskip
\noindent
{\bf {\em Proof.}} Passing to inverses and applying $\varphi$, the relation in Lemma~2.11~a) gives
\[
\begin{split}
&\varphi((1-t_2a_2)^{-1}(1-\rho a_1(t_1)a_2(t_2))^{-1}(1-t_1a_1)^{-1}) \\
&= \frac {h_{a_1}(t_1)h_{a_2}(t_2)}{h_{a_1}(t_1)+h_{a_2}(t_2)-1} h_{a_1+a_2}(t).
\end{split}
\]
On the other hand, taking for instance $b_1 = b_2 = 0$, $s_1 = s_2 = 0$, we are in the bi-free setting in which Lemma~2.9 holds, and we get
\[
\varphi((1-t_2a_2)^{-1}(1-\rho a_1(t_1)a_2(t_2))^{-1}(1-t_1a_1)^{-1}) = h_{a_1}(t_1)h_{a_2}(t_2)
\]
since we will have $b_1(s_1) = b_2(s_2) = 0$ and $h_{b_1}(s_1) = h_{b_2}(s_2) = 1$. This gives then $h_{a_1+a_2}(t) = h_{a_1}(t_1)+h_{a_2}(t_2)-1$. Multiplying this relation with the equality for $t$ in Lemma~2.11~b), we get $th_{a_1+a_2}(t) = t_1h_{a_1}(t_1)$. The remaining equality $t_1h_{a_1}(t_1) = t_2h_{a_2}(t_2)$ is among the assumptions in Lemma~2.11~a).\qed

\bigskip
\noindent
{\bf 2.14. Lemma.} {\em Assume $(a_1,b_1)$ and $(a_2,b_2)$ are bi-free in $({\mathcal A},\varphi)$ and assume $t,s$ are near zero in ${\mathbb C}$. Let further $t_1,t_2,s_1,s_2$ be the numbers near zero provided by Lemma~$2.11$ for $a_1,a_2,t$ and, respectively, $b_1,b_2,s$. Then we have
\[
\frac {h_{a_1}(t_1)h_{b_1}(s_1)}{H_{a_1,b_1}(t_1,s_1)} + \frac {h_{a_2}(t_2)h_{b_2}(s_2)}{H_{a_2,b_2}(t_2,s_2)} - 1 = \frac {h_{a_1+a_2}(t)h_{b_1+b_2}(s)}{H_{a_1+a_2,b_1+b_2}(t,s)}.
\]
}

\bigskip
\noindent
{\bf {\em Proof.}} We shall get the above result by transforming the equality in Lemma~2.9, when the conditions in Lemma~2.11 are satisfied by $t$, $t_1$, $t_2$ and $\rho$ and also the analogous conditions are satisfied by $s$, $s_1$, $s_2$ and $\sigma$ with respect to $b_1,b_2$. The left-hand side of the equality in Lemma~2.9 becomes
\[
\begin{split}
&\frac {h_{a_1}(t_1)h_{a_2}(t_2)}{h_{a_1}(t_1)+h_{a_2}(t_2)-1} \cdot \frac {h_{b_1}(s_1)h_{b_2}(s_2)}{h_{b_1}(s_1)+h_{b_2}(s_2)-1} \cdot H_{a_1+a_2,b_1+b_2}(t,s) \\
&= \frac {h_{a_1}(t_1)h_{a_2}(t_2)h_{b_1}(s_1)h_{b_2}(s_2)}{h_{a_1+a_2}(t)h_{b_1+b_2}(s)} H_{a_1+a_2,b_1+b_2}(t,s)
\end{split}
\]
in view of Lemma~2.13.

On the other hand, the right-hand side in Lemma~2.9 equals
\[
\tfrac {H_{a_1,b_1}(t_1,s_1)H_{a_2,b_2}(t_2,s_2)}{1-(h_{a_1}(t_1)h_{a_2}(t_2)h_{b_1}(s_1)h_{b_2}(s_2))^{-1}(H_{a_1,b_1}(t_1,s_1)-h_{a_1}(t_1)h_{b_1}(s_1))(H_{a_2,b_2}(t_2,s_2)-h_{a_2}(t_2)h_{b_2}(s_2))}.
\]

Dividing by $h_{a_1}(t_1)h_{a_2}(t_2)h_{b_1}(s_1)h_{b_2}(s_2)$ then gives
\[
\begin{split}
&\tfrac {H_{a_1+a_2,b_1+b_2}(t,s)}{h_{a_1+a_2}(t)h_{b_1+b_2}(s)} \\
&= \tfrac {H_{a_1,b_1}(t_1,s_1)H_{a_2,b_2}(t_2,s_2)}{H_{a_1,b_1}(t_1,s_1)h_{a_2}(t_2)h_{b_2}(s_2)+H_{a_2,b_2}(t_2,s_2)h_{a_1}(t_1)h_{b_1}(s_1)-H_{a_1,b_1}(t_1,s_1)H_{a_2,b_2}(t_2,s_2)}.
\end{split}
\]
Taking inverses we get
\[
\frac {h_{a_1+a_2}(t)h_{b_1+b_2}(s)}{H_{a_1+a_2,b_1+b_2}(t,s)} = \frac {h_{a_2}(t_2)h_{b_2}(s_2)}{H_{a_2,b_2}(t_2,s_2)} + \frac {h_{a_1}(t_1)h_{b_1}(s_1)}{H_{a_1,b_1}(t_1,s_1)} - 1
\]
which is the equity we wanted to prove.\qed

\bigskip
With the lemmas we have proved, we will be now able to complete the two steps of the proof of Theorem~2.4, which we outlined in 2.5.

\bigskip
\noindent
{\bf 2.15. Step 1 of the Proof of Theorem 2.4} 

This step will be derived from Lemma~2.14. Remark that $th_a(t) = G_a(t^{-1})$ if $t$ is near $0$, since $t^{-1}$ is then near $\infty$ on the Riemann sphere and $G_a$ is defined near $\infty$. Similarly $ts\,H_{a,b}(t,s) = G_{a,b}(t^{-1},s^{-1})$ if $t$ and $s$ are near $0$.

Next, for $(z,w)$ near $(0,0)$ let $t_j = (K_{a_j}(z))^{-1}$, $s_j = (K_{b_j}(z))^{-1}$, $t = (K_{a_1+a_2}(w))^{-1}$, $s = (K_{b_1+b_2}(w))^{-1}$. Then $t_1,t_2,t$ satisfy the assumptions in Lemma~2.11. Indeed we have $t_jh_{a_j}(t_j) = G_{a_j}(t_j^{-1}) = G_{a_j}(K_j(z)) = z$ for $j = 1,2$ and also $th_{a_1+a_2}(t) = z$. Further, we also have
\[
\begin{aligned}
\frac {t_1h_{a_1}(t_1)}{h_{a_1}(t_1)+h_{a_2}(t_2)-1} &= \frac {z}{\frac {z}{t_1}+\frac {z}{t_2} - 1} \\
&= \frac {1}{K_{a_1}(z)+K_{a_2}(z)-z^{-1}} = \frac {1}{K_{a_1+a_2}(z)} = t.
\end{aligned}
\]
The relations for $s_1,s_2,s$ are checked similarly. With these choices of $t_j,s_j,t,s$, the equality in Lemma~2.14 is an addition property for
\[
\begin{aligned}
\frac {h_{a_j}(t_j)h_{b_j}(s_j)}{H_{a_j,b_j}(t_j,s_j)} - 1 &= \frac {t_jh_{a_j}(t_j)s_jh_{b_j}(s_j)}{t_js_jH_{a_j,b_j}(t_j,s_j)} - 1 \\
&= \frac {G_{a_j}(K_{a_j}(z))G_{b_j}(K_{b_j}(w))}{G_{a_j,b_j}(K_{a_j}(z),K_{b_j}(w))^{-1}} \\
&= \frac {zw}{G_{a_j,b_j}(K_{a_j}(z),K_{b_j}(w))} - 1
\end{aligned}
\]
which is the desired result.

\bigskip
\noindent
{\bf 2.16. Step 2 of the Proof of Theorem 2.4}

It is easy to see that since
\[
\begin{aligned}
(G_{a,b}(z,w))^{-1} &= (z^{-1}w^{-1}H_{a,b}(z^{-1},w^{-1}))^{-1} \\
&= zw\left( \sum_{k \ge 0,l \ge 0} z^{-k}w^{-l}\varphi(a^kb^l)\right)^{-1}
\end{aligned}
\]
the coefficients of $z^pw^q$ in the expansion of $(G_{a,b}(z,w))^{-1}$ near $(\infty,\infty)$ are polynomials in the moments $\varphi(a^kb^l)$. Since $K_a(z) = z^{-1} + R_a(z)$ has also an expansion with coefficients polynomials in $\varphi(a^k)=\varphi(a^kb^0)$ and the same also holds for $K_b(w)$, we infer that the expansion of
\[
zw(G_{a,b}(K_a(z),K_b(w)))^{-1}-1
\]
has as coefficients polynomials in the moments $\varphi(a^kb^l)$.

To prove that the coefficient of $z^kw^l$ is a polynomial of bi-degree $(k,l)$ when $\varphi(a^pb^q)$ is assigned bi-degree $(p,q)$, amounts then to proving that
\[
\begin{split}
&\lambda^{-1}z\mu^{-1}w(G_{\lambda a,\mu b}(K_{\lambda a}(\lambda^{-1}z),K_{\mu b}(\mu^{-1}w)))^{-1} \\
&= zw(G_{a,b}(K_a(z),K_b(w)))^{-1}.
\end{split}
\]
Indeed, we have obviously $h_{\lambda a}(\lambda^{-1}t) = h_a(t)$, $H_{\lambda a,\mu b}(\lambda^{-1}t,\mu^{-1}s) = H_{a,b}(t,s)$, which then give for $G_a(z) = z^{-1}h_a(z^{-1})$ and $G_{a,b}(z,w) = z^{-1}w^{-1}H_{a,b}(z^{-1},w^{-1})$, and that $G_{\lambda a}(\lambda z) = \lambda^{-1}G_a(z)$ and $G_{\lambda a,\mu b}(\lambda z,\mu w) = \lambda^{-1}\mu^{-1}G_{a,b}(z,w)$. Since $K_a$ is the inverse of $G_a$ we have $G_{\lambda a}(\lambda K_a(\lambda z)) = \lambda^{-1}G_a(K_a(\lambda z)) = \lambda^{-1}(\lambda z) = z$ so that $\lambda K_a(\lambda z) = K_{\lambda a}(z)$. This gives $K_{\lambda a}(\lambda^{-1}z) = \lambda K_a(z)$ and hence
\[
\begin{split}
&(G_{\lambda a,\mu b}(K_{\lambda a}(\lambda^{-1}z),K_{\mu b}(\mu^{-1}w)))^{-1} \\
&= (G_{\lambda a,\mu b}(\lambda K_a(z),\mu K_b(w)))^{-1} \\
&= \lambda\mu((G_{a,b}(K_a(z),K_b(w)))^{-1}
\end{split}
\]
which gives the desired result when we multiply by $\lambda^{-1}z\mu^{-1}w$.

We still must show that in the polynomial in moments $\varphi(a^pb^q)$ which is the coefficient of $z^kw^l$, in the expansion of $zw(G_{a,b}(K_a(z),K_b(w)))^{-1} - 1$ the coefficient of the leading term $\varphi(a^kb^l)$ is $1$. We have
\[
\begin{split}
&zw(G_{a,b}(K_a(z),K_b(w)))^{-1} - 1 \\
&= zw((K_a(z))^{-1}(K_b(w))^{-1}H_{a,b}((K_a(z))^{-1},(K_b(w))^{-1}))^{-1}-1.
\end{split}
\]
Let ${\tilde H}_{a,b}(s,t) = 1 - H_{a,b}(s,t) = -\displaystyle{\sum_{\substack{
m \ge 0 \\
n \ge 0 \\
m + n \ge 1}}} s^kt^l\varphi(a^kb^l)$. It is easily seen that in the sense of convergence of formal power series, we have
\[
\begin{split}
&zw(G_{a,b}(K_a(z),K_b(w)))^{-1} - 1 \\
&= -1 + (1+zR_a(z))(1+wR_b(w))\\
&\left( 1 + \sum_{m \ge 0} ({\tilde H}_{a,b}(z(1+R_a(z))^{-1},w(1+R_b(w))^{-1})^m\right).
\end{split}
\]
To check on the coefficients of leading terms, i.e., those terms which involve just one moment, we may leave out in the previous expansion any terms which contain products of $\ge 2$ moments. Thus we arrive at
\[
-1 + (1+zR_a(z)+wR_b(w))(1+{\tilde H}_{a,b}(z(1-R_a(z)),w(1-R_b(w))))
\]
which can then be further reduced to
\[
-1+\left( 1 + \sum_{p \ge 1} z^p\varphi(a^p) + \sum_{q\ge 1} w^q\varphi(b^q)\right)\left( 1 - -\sum_{\substack{
k \ge 0,l \ge 0 \\
k+l > 0}} z^kw^l\varphi(a^kb^l)\right)
\]
and then even further reduced to
\[
\sum_{p\ge 1} z^p\varphi(a^p) + \sum_{q\ge 1} w^q\varphi(b^q) - \sum_{\substack{k\ge 0,l \ge 0 \\
k+l > 0}} z^kw^l\varphi(a^kb^l).
\]
This shows that in order to get the right coefficient of the leading terms, the partial two-variables bi-free $R$-transform is
\[
R_{a,b}(z,w) = 1 - zw(G_{a,b}(K_a(z),K_b(w)))^{-1} + zR_a(z) + wR_b(w).
\]
\qed

\bigskip
\noindent
{\bf 2.17. Remark.} The bi-free cumulants $R_{m,n}(a,b)$ have integer coefficients. This can be seen easily from the fact that the one variable cumulants $R_n(a)$ have integer coefficients and also the expansion we obtained for $R_{a,b}(z,w)$ which involves $G_{a,b}(z,w)$.

\bigskip
\noindent
{\bf 2.18. Remark.} In addition to the relation to the usual one-variable $R$-transform (\cite{8}), note also a certain resemblance in the formula for $R_{a,b}(z,w)$ to formulae for Boolean cumulants (\cite{6}).

\section{Systems with rank $\le 1$ commutation}
\label{sec3}

\noindent
{\bf 3.1. Definition.} An {\em implemented non-commutative probability space} is a triple $({\mathcal A},\varphi,{\mathcal P})$, where $({\mathcal A},\varphi)$ is a non-commutative probability space and ${\mathcal P} = {\mathcal P}^2 \in {\mathcal A}$ is an idempotent so that ${\mathcal P} a{\mathcal P} = \varphi(a){\mathcal P}$. An {\em implemented $C^*$-probability space} $({\mathcal A},\varphi,{\mathcal P})$ will satisfy additional requirements that $({\mathcal A},\varphi)$ be a $C^*$-probability space and that ${\mathcal P} = {\mathcal P}^*$.

\bigskip
\noindent
{\bf 3.2. Example.} If $({\mathcal A},\varphi)$ is a $C^*$-probability space, the GNS construction gives a $*$-representation $\pi: {\mathcal A} \to {\mathcal B}({\mathcal H})$ and a unit vector $\xi \in {\mathcal H}$ so that $\langle \pi(a)\xi,\xi\rangle = \varphi(a)$ if $a \in {\mathcal A}$. Let ${\mathcal P} \in {\mathcal B}({\mathcal H})$ be the orthogonal projector onto ${\mathbb C}\xi$ and let ${\tilde {\mathcal A}}$ be the $C^*$-algebra generated by $\pi({\mathcal A})$ and ${\mathcal P}$ and ${\tilde \varphi}$ the state $\langle \cdot\xi,\xi\rangle$ on ${\tilde {\mathcal A}}$. Then $({\tilde {\mathcal A}},{\tilde \varphi},{\mathcal P})$ is an implemented $C^*$-probability space and $\pi$ gives a morphism of $({\mathcal A},\varphi)$ into $({\tilde {\mathcal A}},{\tilde \varphi})$ so that ${\tilde \varphi}(\pi(a)) = \varphi(a)$. Remark also that the GNS construction can be replaced by a construction in which $\pi$ is a faithful representation, so that $({\mathcal A},\varphi)$ embeds into $({\tilde {\mathcal A}},{\tilde \varphi},{\mathcal P})$.

\bigskip
\noindent
{\bf 3.3. Definition.} A {\em system with rank $\le 1$ commutation} in an implemented non-commutative probability space $({\mathcal A},\varphi,{\mathcal P})$ is a two-faced family $((a_i)_{i \in I},(b_j)_{j \in J})$ in $({\mathcal A},\varphi)$, so that $[a_i,b_j] = \lambda_{ij}{\mathcal P}$ for some $\lambda_{i,j} \in {\mathbb C}$, $i \in I$, $j \in J$. We will call $(\lambda_{ij})_{i \in I,j \in J}$ the coefficients matrix of the system.

\bigskip
\noindent
{\bf 3.4. Remark.} A bipartite system in an implemented non-commutative probability space is the same thing as a system with rank $\le 1$ commutation with coefficients matrix $(\lambda_{ij})_{i \in I,j \in J} = 0$.

\bigskip
\noindent
{\bf 3.5. Proposition.} {\em Let $\mu$ and $\nu$ be the distributions of two systems with rank $\le 1$ commutation, with index sets $(I,J)$ and, respectively, $(K,L)$ for the variables.}

a) {\em Let $\pi$ be the distribution of a two-faced system of variables indexed by $(I \coprod K,J \coprod L)$ so that if $\pi',\pi''$ are the restrictions of $\pi$ to the index pairs $(I,J)$ and $(K,L)$, then $\pi' = \mu$, $\pi'' = \nu$ and $\pi',\pi''$ are bi-free when seen as parts of the distribution $\pi$. Then $\pi$ is the distribution of a system with rank $\le 1$ commutation and coefficients matrix $(\lambda_{ij})_{i \in I,j \in J} \oplus (\lambda_{kl})_{k \in K,l \in L}$ the direct sum of the coefficients matrices for the systems with distributions $\mu$ and $\nu$.}

b) {\em If $I = K$, $J = L$ then the bi-free convolution $\mu \boxplus\boxplus \nu$ is the distribution of a system with rank $\le 1$ commutation and coefficients matrix the sum of the coefficients matrices for the systems with distributions $\mu$ and $\nu$.
}

\bigskip
\noindent
{\bf {\em Proof.}} a) Let $((a'_i)_{i \in I},(b'_j)_{j \in J})$ in $({\mathcal A}',\varphi',{\mathcal P}')$ and $((a''_k)_{k \in K},(b''_l)_{l \in L})$ in $({\mathcal A}'',\varphi'',{\mathcal P}'')$ be the systems with distributions $\mu,\nu$ and let $(\lambda'_{ij})_{i \in I,j \in J}$, $(\lambda''_{kl})_{k \in K,l \in L}$ be the corresponding coefficients matrices (a slight deviation from the notation without primes and double primes in the statement of proposition). Let further $({\mathcal X}',\overset{\circ}{{\mathcal X}}{}',\xi')$ and $({\mathcal X}'',\overset{\circ}{{\mathcal X}}{}'',\xi'')$ denote the vector spaces with specified state-vectors $({\mathcal A}'{\mathcal P}',(I-{\mathcal P}'){\mathcal A}'{\mathcal P}',{\mathcal P}')$ and $({\mathcal A}''{\mathcal P}'',(I-{\mathcal P}''){\mathcal A}''{\mathcal P}'',{\mathcal P}'')$ endowed with the natural left actions of ${\mathcal A}'$ and ${\mathcal A}''$. Consider $({\mathcal X},\overset{\circ}{{\mathcal X}},\xi) = ({\mathcal X}',\overset{\circ}{{\mathcal X}}{}',\xi') * ({\mathcal X}'',\overset{\circ}{{\mathcal X}}{}'',\xi'')$ and $\lambda',\rho',\lambda'',\rho''$ the left and right actions of ${\mathcal L}({\mathcal X}')$ and ${\mathcal L}({\mathcal X}'')$ on ${\mathcal X}$ viewed as homomorphisms into ${\mathcal L}({\mathcal X})$. The actions of ${\mathcal A}'$ and ${\mathcal A}''$ on ${\mathcal X}'$ and ${\mathcal X}''$, respectively, yield homomorphisms $\eta': {\mathcal A}' \to {\mathcal L}({\mathcal X}')$, $\eta'': {\mathcal A}'' \to {\mathcal L}({\mathcal X}'')$. It is easily seen that denoting by ${\mathcal P}$ the idempotent in ${\mathcal L}({\mathcal X})$ so that ${\mathcal P}{\mathcal X} = {\mathbb C}\xi$ and $\ker {\mathcal P} = \overset{\circ}{{\mathcal X}}$, we have
\[
\begin{aligned}
[\lambda'(\eta'(a'_i)),\rho'(\eta'(b'_j))] &= \lambda'_{ij} {\mathcal P}, \\
[\lambda''(\eta''(a''_k)),\rho''(\eta''(b''_l))] &= \lambda''_{kl} {\mathcal P}, \\
[\lambda'(\eta'(a'_i)),\rho''(\eta''(b''_l))] &= 0 \\
[\lambda''(\eta''(a''_k)),\rho'(\eta'(b'_j))] &= 0.
\end{aligned}
\]
Then $((\lambda'(\eta'(a'_i)))_{i \in I} \coprod (\lambda''(\eta''(a''_k)))_{k \in K}$, $(\lambda'(\eta'(b'_j)))_{j \in J} \coprod (\lambda''(\eta''(b''_l)))_{l \in L})$ is a system with rank $\le 1$ commutation in $({\mathcal L}({\mathcal X}),\varphi_{\xi},{\mathcal P})$, which from the distributions point of view is made of bi-free copies of $\mu$ and $\nu$ and has coefficients matrix $(\lambda'_{ij})_{i \in I,j \in J} \oplus (\lambda''_{kl})_{k \in K,l \in L}$.

b) We can use the construction in the proof of a). Then the system $((\lambda'(\eta'(a'_i)) + \lambda''(\eta''(a''_i)))_{i \in I}$, $(\rho'(\eta'(b'_j)) + \rho''(\eta''(b''_j)))_{j \in J})$ in $({\mathcal L}({\mathcal X}),\varphi_{\xi},{\mathcal P})$ is a system with rank $\le 1$ commutation with distribution $\mu \boxplus\boxplus \nu$. Since
\[
[\lambda'(\eta'(a'_i)) + \lambda''(\eta''(a''_i)),\ \rho'(\eta'(b'_j)) + \rho''(\eta''(b''_j))] = \lambda'_{ij} + \lambda''_j
\]
the coefficients matrix is clearly $(\lambda'_{ij} + \lambda''_{ij})_{i \in I,j \in J}$.\qed

\bigskip
\noindent
{\bf 3.6.} We will introduce in this subsection the terminology of {\em bands in index sequences}. For systems with rank $\le 1$ commutations, we will have to deal with up to two bands, but there is no extra cost to give general definitions, which may provide a better perspective.

If $(I,J)$ is a pair of index sets for a two-faced system of non-commutative random variables, a map $\alpha: \{1,\dots,m\} \to I \coprod J$ will be called an {\em index sequence of length $m$}. A {\em band in the index sequence $\alpha$} is a maximal interval $[p,q]$, $1 \le p \le q \le m$, $p,q \in {\mathbb N}$, so that all $\alpha(r)$, $p \le r \le q$ are in only one of the two sets $I$ and $J$, that is the corresponding variables will be either all left or all right variables. An index sequence $\alpha$ which has $n$ bands will also be called an {\em $n$-band index sequence}. It may also be useful to say $\alpha$ {\em starts left} if $\alpha(1) \in I$ and, hence, also its first band is left or $\alpha$ {\em starts right} if $\alpha(1) \in J$. Similarly we may say $\alpha$ {\em ends left} or $\alpha$ {\em ends right} and clearly $\alpha$ starts and ends on the same side if the band number is odd and starts and ends on different sides if the and number is even. If $\alpha$ has $n$ bands the corresponding non-commutative moment will be called {\em $n$ bands moments} and the corresponding cumulants will also be called {\em $n$ bands cumulants}. Also the part of the distribution of a two-faced system of non-commutative random variables involving moments with $\le n$ bands, will be called its {\em $n$ bands distribution}. Clearly the $n$ bands distribution can also be given by the cumulants with $\le n$ bands. In particular, {\em the bi-free additive convolution gives rise to operations on the corresponding $n$ bands distributions}. It may also be useful at times to consider a subset of the moments or cumulants with $\le n$ bands by requiring that those with precisely $n$ bands start left or right, etc. We may thus speak of the {\em $n$ bands distribution starting left} (we repeat that the starting condition applies only to the moments or cumulants with maximum number of bands, that is with $n$ bands).

\bigskip
\noindent
{\bf 3.7. Proposition.} {\em The distribution of a system with rank $\le 1$ commutation is completely determined by its $2$-bands distribution starting left and its coefficients matrix.
}

\bigskip
\noindent
{\bf {\em Proof.}} Let $({\mathcal A},\varphi,{\mathcal P})$ be an implemented non-commutative probability space and let $((z_i)_{i \in I},(z_j)_{j \in J})$ be a system with rank $\le 1$ commutation in $({\mathcal A},\varphi,{\mathcal P})$. Let further $IJ = \{(i_1,\dots,i_p,j_1,\dots,j_q) \mid p \ge 0,q \ge 0,i_1,\dots,i_p \in I,j_1,\dots,j_q \in J\}$ which contains a special ``empty'' element with $p = q = 0$. If $k \in I \coprod J$ we will denote by $T_k$ the linear operator of right multiplication by $z_k$ on ${\mathcal A}$, that is $T_ka = az_k$ if $a \in {\mathcal A}$. By ${\mathcal V} \subset {\mathcal A}$ we will denote the subspace spanned by the ${\mathcal P} z_{i_1}\dots z_{i_p}z_{j_1}\dots z_{j_q}$ where $(i_1,\dots,i_p,j_1,\dots,j_q) \in IJ$ which means that we assume ${\mathcal P} \in {\mathcal V}$ corresponding to the empty element of $IJ$. We claim that $T_k{\mathcal V} \subset {\mathcal V}$. If $k \in J$ this is obvious since $(i_1,\dots,i_p,j_1,\dots,j_q,k) \in IJ$ in this case and hence $T_k{\mathcal P} z_{i_1}\dots z_{i_p}z_{j_1}\dots z_{j_q} = {\mathcal P} z_{i_1}\dots z_{i_p}z_{j_1}\dots z_{jq}z_k \in {\mathcal V}$. If $k \in I$ we have
\[
\begin{aligned}
T_k{\mathcal P} z_{i_1}\dots z_{i_p}z_{j_1}\dots z_{j_q} &= {\mathcal P} z_{i_1}\dots z_{i_p}z_kz_{j_1}\dots z_{j_q} \\
&+ \sum_{t=1}^q {\mathcal P} z_{i_1}\dots z_{i_p}z_{j_1}\dots z_{j_{t-1}}[z_{j_t},z_k]z_{j_{t+1}}\dots z_{j_q} \\
&= {\mathcal P} z_{i_1}\dots z_{i_p}z_kz_{j_1}\dots z_{j_q} \\
&- \sum_{t=1}^q \varphi(z_{i_1}\dots z_{i_p}z_{j_1}\dots z_{j_{t-1}})\lambda_{k,j_t} {\mathcal P} z_{j_{t+1}}\dots z_{j_q}
\end{aligned}
\]
since $[z_{j_t},z_k] = -\lambda_{k,j_t}{\mathcal P}$ and 
\[
{\mathcal P} z_{i_1}\dots z_{i_p}z_{j_1}\dots z_{j_{t-1}} {\mathcal P} = \varphi(z_{i_1}\dots z_{i_p}z_{j_1}\dots z_{j \in 1}){\mathcal P}.
\]
This shows that $T_k{\mathcal V} \subset {\mathcal V}$ and, moreover, the formulae for computing $T_k{\mathcal P} z_{i_1}\dots z_{i_p}z_{j_1}\dots z_{j_q}$ are determined by the $2$-bands distribution starting left and the coefficients matrix of the system with rank $\le 1$ commutation. The proposition then follows by applying several times this fact to compute ${\mathcal P} z_{k_1}\dots z_{k_r} = T_{k_r} \dots T_{k_1}{\mathcal P}$ where $k_1,\dots,k_r \in I \coprod J$ and then passing to $\varphi(z_{k_1}\dots z_{k_r})$ which is given by
\[
\varphi(z_{k_1}\dots z_{k_r}){\mathcal P} = (T_{k_r}\dots T_{k_1}{\mathcal P}){\mathcal P}
\]
and since $T_{k_1}\dots T_{k_r}{\mathcal P} \in {\mathcal V}$ the computation of $(T_{k_1}\dots T_{k_r}{\mathcal P}){\mathcal P}$ gives a polynomial in
\[
{\mathcal P} z_{i_1}\dots z_{i_p}z_{j_1}\dots z_{j_q}{\mathcal P} = \varphi(z_{i_1}\dots z_{i_p}z_{j_1}\dots z_{j_q}){\mathcal P}.
\]
Clearly at each step the computations involve only the $2$-bands distribution starting left and the coefficients matrix, which proves the proposition.\qed

\bigskip
\noindent
{\bf 3.8. Example.} The {\em dual systems} in section~5 of \cite{9} provide examples of systems with rank $\le 1$ commutation. Given $({\mathcal M},\tau)$ a von~Neumann algebra with a normal faithful trace-state, $I \in {\mathcal B} \subset {\mathcal M}$ a $*$-subalgebra and $X_j = X_j^* \in {\mathcal M}$, $1 \le j \le n$, a dual system to $(X_1,\dots,X_n;{\mathcal B})$ in $L^2({\mathcal M},\tau)$ is an $n$-tuple $(Y_1,\dots,Y_n)$ of operators $Y_j = Y_j^* \in {\mathcal B}(L^2({\mathcal M},\tau))$, $1 \le j \le n$ such that $[{\mathcal B},Y_j] = 0$ and $[X_j,Y_k] = i \delta_{jk} {\mathcal P}$ where ${\mathcal P}$ is the orthogonal projection onto ${\mathbb C}1$. Clearly $((X_1,\dots,X_n) \coprod (b_k)_{k \in K})$, $(Y_1,\dots,Y_n)$ where $(b_k)_{k \in K}$ is any family of elements in ${\mathcal B}$, is a system with rank $\le 1$ commutation in $({\mathcal B}(L^2({\mathcal M},\tau)),\omega,{\mathcal P})$ where $\omega$ is the state $\langle \cdot 1,1\rangle$ on ${\mathcal B}(L^2({\mathcal M},\tau))$. The coefficient matrix in this case is $\lambda_{pq} = i\delta_{pq}$ if $1 \le p,q \le n$ and $\lambda_{pk} = 0$ if $1 \le p \le n$, $k \in K$. There are additional features in this example, such as the framework and the distribution being $C^*$, the $X_p$'s and $Y_q$'s being hermitian and the joint distribution of $X_p$'s and $(b_k)_{k \in K}$ being a tracial $*$-distribution. Proposition~3.5 can be viewed as an algebraic version and generalization of results for dual systems in \cite{9}.

\bigskip
\noindent
{\bf 3.9. Example.} Let ${\mathcal H}$ be a Hilbert space and let $T \in {\mathcal B}({\mathcal H})$ be a hyponormal operator with rank one self-commutator, that is $T^*T - TT^* = \lambda {\mathcal P}$, $\lambda > 0$ where ${\mathcal P}$ is the orthogonal projection onto an one-dimensional subspace ${\mathbb C}\xi \in {\mathcal H}$, $\|\xi\| = 1$. In $({\mathcal B}({\mathcal H}),\omega,{\mathcal P})$ where $\omega(S) = \langle S\xi,\xi\rangle$ for $S \in {\mathcal B}({\mathcal H})$, one can derive from $T$ examples of pairs with rank~1 self-commutator. Indeed if $\Omega = \begin{pmatrix}
a & b \\
c & d
\end{pmatrix} \in GL(2;{\mathbb C})$ we can consider $(aT + bT^*,cT + dT^*)$ which has commutator
\[
[aT+bT^*,cT+dT^*] = -\det(\Omega)\lambda{\mathcal P}.
\]
In case $\Omega = \begin{pmatrix}
1 & 1 \\
-i & i
\end{pmatrix}$ we get a pair of hermitian operators, like in the simplest case of the dual systems in 3.8. Clearly, the bi-free additive convolutions of the distributions of such pairs are distributions of the same kind, that is arising from a hyponormal operator. On the other hand, taking for instance $\Omega = \begin{pmatrix}
1 & 0 \\
0 & 1
\end{pmatrix}$ we get the pair $(T,T^*)$ and bi-free addition of such pairs leads to a pair which does not correspond in the same straightforward way to a hyponormal operator with rank one self-commutator.

\bigskip
\noindent
{\bf 3.10. Example.} The bi-free Gaussian two-faced system in Theorem~7.4 of Part~I is also an example with rank $\le 1$ commutation. Indeed, if ${\mathcal H}$ is a Hilbert space, ${\mathcal T}({\mathcal H})$ the full Fock space and $h,h^*: I \coprod J \to {\mathcal H}$ then $z_i = l(h(i)) + l^*(h^*(i))$, $i \in I$ and $z_j = r(h(j)) + r^*(h^*(j))$, $j \in J$ has rank $\le 1$ commutation in $({\mathcal B}({\mathcal T}({\mathcal H})),\omega,{\mathcal P})$ where $\omega(T) = \langle T1,1\rangle$, ${\mathcal P}\eta = \langle \eta,1\rangle 1$. We have $[z_i,z_j] = (\langle h(j),h^*(i)\rangle - \langle h(i),h^*(j)\rangle){\mathcal P}$.

\section{Classically independent left and right faces}
\label{sec4}

\noindent
{\bf 4.1. Theorem.} {\em Let $\pi = ((a_i)_{i \in I},(b_j)_{j \in J})$ be a two-faced system of non-commutative random variables in $({\mathcal A},\varphi)$.
}

a) {\em The distribution of $\pi$ is equal to the distribution of $\pi' = ((a'_i)_{i \in I},(b'_j)_{j \in J})$ where $a'_i = a_i$, $b'_j = 0$ iff the only non-vanishing bi-free cumulants of $\pi$ correspond to index sequences $\alpha$ the range of which is in $I$. The obvious analogue for index sequences with range in $J$ also holds.
}

b) {\em The distribution of $\pi$ has the property that the only non-vanishing bi-free cumulants correspond to index sequences with ranges either in $I$ or in $J$ iff the families of non-commutative random variables $(a_i)_{i \in I}$ and $(b_j)_{j \in J}$ are classically independent in $({\mathcal A},\varphi)$.
}

\bigskip
\noindent
{\bf {\em Proof.}} a) Proposition~5.6 of Part~I implies a multi-homogeneity indexed by $I \coprod J$ for bi-free addition, which implies that the bi-free cumulants in view of their uniqueness are also multi-homogeneous. More precisely, if $\alpha: \{1,\dots,n\} \to I \coprod J$, for every $k \in I \coprod J$ assign the degree $\#\alpha^{-1}(k)$ and for every moment of $\pi$ assign the corresponding multi degree $I \coprod J \ni k \to \#\alpha^{-1}(k) \in {\mathbb Z}_{\le 0}$ where $\alpha$ is the index sequence of the monomial the expectation of which is the moment. Then the bi-free cumulant is multi-homogeneous with degrees map corresponding to its index sequence. In view of this, $\pi'$ is obtained from $\pi$ by multiplying the $b_i$'s by $0$, we infer that the corresponding cumulants are obtained from those of $\pi$ by multiplying the ones involving right indices by $0$. This proves the only if part and also the if part since the cumulants completely determine the distribution.

b) Like in a) consider besides $\pi'$ also $\pi''$ where $a''_i = 0$ and $b''_j = b_j$. In view of Proposition~2.16 of Part~I the bi-free sum of $\pi'$ and $\pi''$ consists of classically independent copies of $(a_i)_{i \in I}$ and $(b_j)_{j \in J}$ and in view of the addition property of cumulants its bi-free cumulants are either zero or free cumulants of $(a_i)_{i \in I}$ or $(b_j)_{j \in J}$. This proves the if part of b) and the only if actually also follows since cumulants completely determine the distribution.\qed

\bigskip
\noindent
{\bf 4.2. Corollary.} {\em If $\mu$ and $\nu$ are distributions for two-faced families with classically independent left and right variables and same index sets, then $\mu \boxplus\boxplus \nu$ is also a distribution with classically independent left and right variables.
}


\begin{thebibliography}{99}

\bibitem{1} Clancey, K. {\em Seminormal operators}, Lecture Notes in Math. {\bf 742}, Springer (1979).

\bibitem{2} Haagerup, U., {\em On Voiculescu's $R$- and $S$-transforms for free non-commuting variables}, in ``Free Probability Theory'', D.-V. Voiculescu, editor, Fields Institute Communications {\bf 12}, pp.~127--148, AMS (1997).

\bibitem{3} Martin, M., and Putinar, M., {\em Lectures on hyponormal operators}, Operator Theory: Advances and Applications {\bf 39}, Birkh\"auser (1989).

\bibitem{4} Mineev-Weinstein, M.; Putinar, M.; and Teodorescu, R., {\em Random matrices in $2D$ Laplacian growth and operator theory}, J. Phys. A {\bf 41} (2008), No.~26, 263001, 74~pp.

\bibitem{5} Speicher, R., {\em Combinatorial theory of the free product and operator-valued free probability theory}, Memoirs of the AMS {\bf 132} (1998), No.~627.

\bibitem{6} Speicher, R., and Woroudi, R., {\em Boolean convolution}, in ``Free Probability Theory'', D.-V. Voiculescu, editor, Fields Institute Communications {\bf 12}, pp.~267--279, AMS (1997).

\bibitem{7} Voiculescu, D.-V., {\em Symmetries of some reduced free product $C^*$-algebras}, in ``Operator Algebras and their Connections with Topology and Ergodic Theory'', Lecture Notes in Math. {\bf 1132} (1985), Springer Verlag, 556--588.

\bibitem{8} Voiculescu, D.-V., {\em Addition of certain non-commuting random variables}, J. Funct. Anal. {\bf 66} (1986), 323--346.

\bibitem{9} Voiculescu, D.-V., {\em The analogues of entropy and of Fisher's information measure in free probability theory V. Non-commutative Hilbert transforms}, Inventiones Math. {\bf 132} (1998), 189--227.

\bibitem{10} Voiculescu, D.-V., {\em Free probability for pairs of faces I}, preprint (2013), arXiv: 1306.6082, 33~pp.

\bibitem{11} Voiculescu, D.-V.; Dykema, K.~J.; and Nica, A., {\em Free Random Variables}, CRM Monograph Series {\bf 1}, AMS (1992). 

\end{thebibliography}
\end{document}